\documentclass[10p, twocolumn]{article}
\usepackage{amsfonts, amsthm, amsmath}
\usepackage{graphicx}
\usepackage{amssymb}
\usepackage{epstopdf}

\newtheorem{theorem}{Theorem}

\newtheorem{Def}{Definition}

\newtheorem{Q}{Question}

\newcommand{\Realline}{\mathbb{R}}
\newcommand{\Integers}{\mathbb{Z}}

\newcommand{\CC}{\mathbb{C}}
\newcommand{\Foliation}{\mathcal{F}}
\newcommand{\Circle}{\mathrm{S}^1}

\title{Illumination problems on translation surfaces with planar infinities}

\date{}

\author{Nikolay Dimitrov \thanks{Department of Mathematics and Statistics,
        McGill University, {\tt dimitrov@math.mcgill.ca}}}

\index{Dimitrov, Nikolay}

\begin{document}
\thispagestyle{empty}

\maketitle

\begin{abstract}
In the current article we discuss an illumination problem proposed
by Urrutia and Zaks. The focus is on configurations of finitely
many two-sided mirrors in the plane together with a source of
light placed at an arbitrary point. In this setting, we study the
regions unilluminated by the source. In the case of rational-$\pi$
angles between the mirrors, a planar configuration gives rise to a
surface with a translation structure and a number of planar
infinities. We show that on a surface of this type with at least
two infinities, one can find plenty of unilluminated regions
isometric to unbounded planar sectors. In addition, we establish
that the non-bijectivity of a certain circle map implies the
existence of unbounded dark sectors for rational planar mirror
configurations illuminated by a light-source.
\end{abstract}

\section{Introduction} Consider a planar domain with a light
reflecting boundary. Place a source of light at a point inside the
domain. Assume that the source emits rays in all directions. Each
ray follows a straight line and whenever it reaches the boundary
it is reflected according to the rule that the angle of incidence
equals the angle of reflection. A point from the domain is
considered \emph{illuminated} by the source whenever there is a
ray that reaches the point either directly or after a series of
reflections. In this setting, one can ask the following questions,
also known as \emph{illumination problems}.

\begin{Q} If we place the source of light at any point in the domain, will
all of the domain be illuminated? If not, what could be said about
the non-illuminated regions? \end{Q}

\begin{Q} Is there a point from which the light source can illuminate
the entire domain?\end{Q}

These problems are often attributed to E. Straus who posed them
sometime in the early fifties and first published by V. Klee in
1969 \cite{T}. Some famous examples and interesting results are
Penrose's room \cite{CFG}, Tokarsky's example \cite{T} as well as
the article \cite{HST} by Hubert, Schmoll and Troubetzkoy on
illumination on Veech surfaces.

In 1991, J. Urrita and J. Zaks proposed the following problem
\cite{U}. Assume we are given a finite number of disjoint compact
line segments in the plane each representing a mirror that
reflects light on both sides (a two-sided mirror).
\begin{figure}[h!]
\centering
\includegraphics[width=7.9cm]{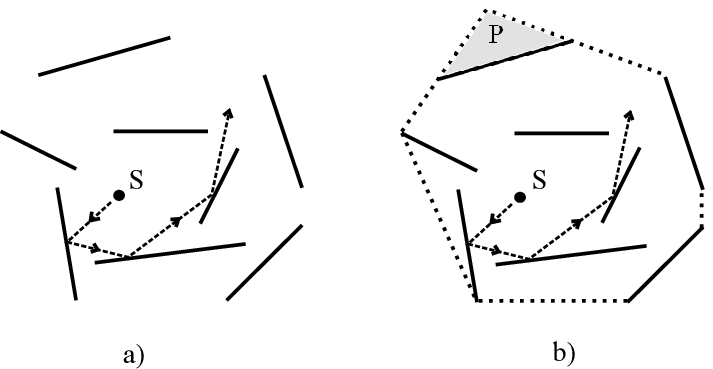}
\caption{} \label{Mirrors}
\end{figure}
Let $p_0$ be
any point on the plane not incident to any of the segments. Then,
the complement of the set of mirrors is an unbounded domain with
light-reflecting boundary and if we place a source of light $S$ at
$p_0$ we can pose questions 1 and 2. Figure \ref{Mirrors}a depicts
an example of a two-sided mirror configuration with a light
emitting source $S$. The convex hull of the mirrors is a polygon.
If $S$ is in the convex hull, one can construct a triangle $P$
unilluminated by $S$, like the shaded one on figure
\ref{Mirrors}b. To do that, it is sufficient for a mirror segment
to be an edge of the convex hull.

In this paper we are interested in finite two-sided mirror
configurations with the following property: any pair of lines
determined by the mirror segments are either parallel or intersect
at an angle which is a rational multiple of $\pi.$ We will call
such a configuration a \emph{rational mirror configuration} and
the domain obtained as a complement of the mirrors will be called
\emph{rational mirror domain}. For those, we will find conditions
that will guarantee the existence of unbounded unilluminated
sectors in the plane (see definition \ref{Definition sector}).

A rational mirror domain can be "unfolded" into a surface that
carries a flat structure with conical singularities and trivial
holonomy group \cite{HS}, \cite{M}. This means that the surface
has a special atlas with the property that away from the cone
points, the transition maps between two charts from the atlas are
Euclidean translations. In the literature, such an object is
called \emph{a translation surface}. As a result, the piecewise
linear trajectory of a light ray in the original domain becomes a
smooth geodesic on the flat surface. Thus, one can think of a
light source placed at a nonsingular point on the surface,
emitting geodesic rays in all directions. Any other point on the
surface is considered \emph{illuminated} if there is a smooth
geodesic connecting the source to the point. In this way, one can
ask questions 1 and 2 for the translation surface. Notice that
there are regions on it isometric to complements of compact sets
in the plane. We will call a surface with such geometry \emph{a
translation surface with planar infinities}.

A translation surface with planar infinities gives rise to a pair
$(X,\omega)$ where $X$ is a closed surface with a complex
structure and $\omega$ is a meromorphic differential on $X$ with
only double poles and zero residues. The zeroes of $\omega$ are
the cone points of the flat structure \cite{HS},\cite{M}, and
around each pole the surface looks like the complement of a
compact set in the plane. The converse is also true. A pair
$(X,\omega)$ of a closed Riemann surface and a meromorphic
differential with only double poles and zero residues induces a
translation structure on $X$ with planar infinities.

\begin{Def} \label{Definition_surface}
The pair $(X,\omega)$ is called a translation surface with planar
infinities whenever the following conditions hold:

\smallskip
\noindent \emph{(1)} $X$ is a closed surface with a complex
structure;

\smallskip
\noindent \emph{(2)} $\omega$ is a meromorphic differential on
$X$;

\smallskip
\noindent \emph{(3)} Every pole of $\omega$ is of order exactly
$2$ and the residue at that pole is zero. We will refer to the
poles of $\omega$ as planar infinities.
\end{Def}

In this study we will be interested in a special type of domains
both on a translation surface with planar infinities and in the
plane.

\begin{figure}[h!] \centering
\includegraphics[width=8cm]{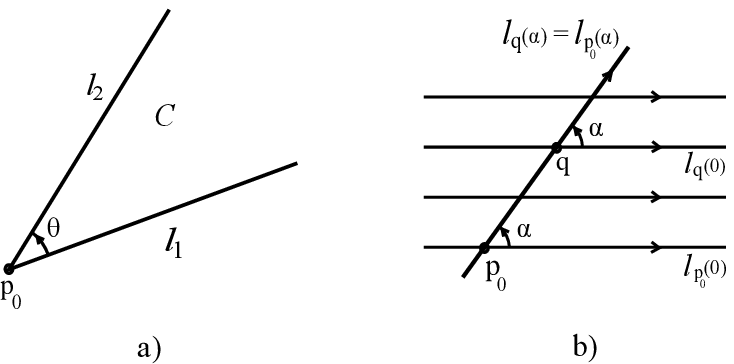}
\caption{} \label{Sector}
\end{figure}

\begin{Def} \label{Definition sector}
\emph{a)} Let $l_1$ and $l_2$ be two half-lines in the plane both
starting form a point $p_0$ and going to infinity. Let $\theta$ be
the angle between $l_1$ and $l_2$ at  the vertex $p_0$, measured
counterclockwise from $l_1$ to $l_2$. Then, the open region $C$
bounded by $l_1$ and $l_2$, whose internal angle at $p_0$ is
$\theta$, is called an infinite sector of angle $\theta$ (see
figure \ref{Sector}a).

\smallskip
\noindent \emph{b)} An open subdomain $C$ of a translation surface
with planar infinities $(X,\omega)$ is called an infinite sector
of angle $\theta$ whenever there exists a chart from the
translation atlas of $(X,\omega)$ that maps $C$ isometrically to a
planar infinite sector of angle $\theta$ like the one defined in
point a.
\end{Def}
On any translation surface $(X,\omega)$ one can always find an
orientable foliation $\Foliation_{\omega}$ with singularities,
whose leaves are geodesics. Indeed, let us foliate the Euclidean
plane into horizontal straight lines, oriented as usual from left
to right. Since each transition map between two charts is a
Euclidean translation, it sends horizontal lines to horizontal
lines (line orientation preserved). Thus, pulling back onto the
surface the planar horizontal foliation from all translation
charts defines globally the desired foliation
$\Foliation_{\omega}$. Moreover, the singularities of
$\Foliation_{\omega}$ are the cone points of the surface
$(X,\omega)$, i.e. the zeroes of the differential $\omega.$ We
call $\Foliation_{\omega}$ the \emph{horizontal foliation} of the
surface and its leaves - the \emph{horizontal geodesics} of the
surface. At each non-singular point $p_0$ of $(X,\omega)$ the
oriented horizontal geodesic $l_{p_0}(0)$ from
$\Foliation_{\omega}$ defines a \emph{positive horizontal
direction} at $p_0$. The counterclockwise angle $\alpha$ between
$l_{p_0}(0)$ and an arbitrary oriented geodesic $l_{p_0}(\alpha)$
through $p_0$ is called the \emph{direction} of $l_{p_0}(\alpha)$
at $p_0$ (see figure \ref{Sector}b). From now on,
$l_{p_0}(\alpha)$ denotes the geodesic ray on $(X,\omega)$
starting from $p_0 \in X$ and going in the direction of angle
$\alpha$. It is important to emphasize that, since we are working
with a translation surface, the intersection of the geodesic
$l_{p_0}(\alpha)$ with any other horizontal geodesic $l_q(0)$ will
always form the same angle $\alpha,$ as shown locally on figure
\ref{Sector}b. In other words, just like in the plane, a geodesic
on $(X,\omega)$ does not changes its angle with respect to the
horizontal direction. Since a direction at any non-singular point
$p \in X$ is defined as an angle $\alpha \in \Realline \mod 2
\pi$, we can identify the set of all directions at $p$ with the
unit circle $\Circle = \{z \in \CC \,:\,\,|z|=1\}$. The point $1
\in \Circle$ gives the horizontal direction $\alpha=0$.

\section{Results} It is natural to ask questions about the
behavior of the geodesics on a surface. The first question we will
address is the following. \emph{On a translation surface with
planar infinities, where do most geodesics emanating from a
nonsingular point go?} As it turns out, almost all of them fall
onto the poles of the surface. Same is true for any rational
mirror configuration in the plane.

\begin{theorem} \label{Main theorem about directions} The following
two statements are true:

\smallskip
\noindent \emph{(1)} Let $(X,\omega)$ be a translation surface
with planar infinities and let $p_0 \in X$ be non-singular. Then
the set of all directions $\alpha \in \Circle$ for which the
geodesic passing through $p_0$ in direction of $\alpha$ goes to
one of the poles of $\omega$ is open and dense in the circle
$\Circle$;

\smallskip
\noindent \emph{(2)} Assume we are given a rational mirror
configuration in the plane and let $p_0$ be a point not lying on
any of the mirrors. Then the set of all directions $\alpha \in
\Circle$ for which the piece-wise linear reflected trajectory
starting from $p_0$ in the direction of $\alpha$ goes to infinity
is open and dense in the circle $\Circle$.
\end{theorem}

The next result establishes the existence of infinite
unilluminated sectors and large unbounded regions on translation
surfaces with more than one planar infinity.

\begin{theorem} \label{Main theorem illumination on surface}
Let $(X,\omega)$ be a translation surface with at least two planar
infinities. Then, for any point $p_0$ on $X \setminus
(zeroes(\omega) \cup \text{poles}(\omega))$ there exists an
infinite sector $C$ on $(X,\omega)$ unilluminated by $p_0$, i.e.
for any point $p \in C$ there is no smooth geodesic on
$(X,\omega)$ that connects $p_0$ to $p$. Moreover, there exists a
region on $(X,\omega)$ consisting of unilluminated,
non-overlapping infinite sectors of total angle $2 \pi (k-1)$,
where $k$ is the number of poles of $\omega$.
\end{theorem}

The main ideas used in the proof of theorem \ref{Main theorem
illumination on surface} can be adjusted to the study of
illumination problems for rational mirror configurations in the
plane. For instance, an interesting question put in an every day
language, is the following. \emph{How big of an object can be
hidden from a stationary observer in a rational mirror domain? Can
we hide a car? A whole parking lot of cars?} Precisely speaking,
we would like to find a basic condition that will ensure the
existence of an infinite unilluminated sector for a light source
placed at a point inside a rational mirror domain.

Let $D$ be a rational mirror domain and let $p_0 \in D$. Draw a
large enough circle $K$, so that its interior contains the mirrors
from the configuration and the light source at the point $p_0$.
Denote by $U_{p_0}$ the open dense set of all directions which go
to infinity, provided by theorem \ref{Main theorem about
directions}. For an angle $\alpha \in U_{p_0} \subset \Circle$
follow the straight line $l_{p_0}(\alpha)$ starting form $p_0$ in
direction of $\alpha$. Whenever the line reaches a mirror it is
reflected, changing its direction. In this way, a piecewise linear
trajectory is formed,\begin{figure} \centering
\includegraphics[width=7.9cm]{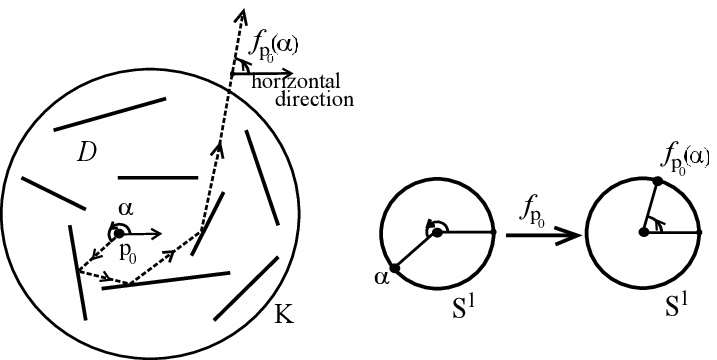}
\caption{} \label{PictureCircleMap}
\end{figure} which at some point leaves the disc bounded
by $K$ never to come back to it. Denote by $f_{p_0}(\alpha)$ the
angle between the horizontal direction of $\CC$ and the portion of
the trajectory that is outside the circle $K$. As a result, we
obtain a map $f_{p_0} : U_{p_0} \longrightarrow \Circle$. For a
picture of the construction of $f_{p_0}$ see figure
\ref{PictureCircleMap}.

The map $f_{p_0}$ is defined almost everywhere on the unit circle.
In fact, its domain $U_{p_0}$ is open and dense in $\Circle$.
Moreover, $f_{p_0}$ is a rotation when restricted to any connected
component of $U_{p_0}$. Our hope is that finding ways to study the
combinatorial properties of $f_{p_0}$ may facilitate the search
for unbounded unilluminated sectors in rational mirror domains.

\begin{theorem} \label{Main theorem illumination in plane}
Assume we are given a rational mirror configuration. For an
arbitrary point $p_0$, not located on any of the mirrors,
construct the circle map $f_{p_0}$ as explain in the previous two
paragraphs (see also figure \ref{PictureCircleMap}). If $f_{p_0}$
is not injective, then there exists an infinite sector in the
plane unilluminated by $p_0$.
\end{theorem}

\section{Translation surfaces.}

In the current section we discuss translation surfaces and show
how to construct one from a rational mirror configuration. To
illustrate the idea better, we apply the procedure to an example.

\paragraph{Various descriptions.} A translation surface is a
closed surface $X$ with a finite set of points $\Sigma \subset X$,
called singularities, and a cover of $X \setminus \Sigma$ by open
charts $$\{(W_{a},\varphi_a)\,\,| \,\, W_a \subseteq X \setminus
\Sigma\, , \,\, \varphi_a : W_a \to \CC\}$$ having the property
that whenever $W_a \cap W_b \neq \varnothing$ the transition map
between the two charts $(W_a,\varphi_a)$ and $(W_b,\varphi_b)$ is
a Euclidean translation, i.e. $z_b = \varphi_b^{-1} \circ
\varphi_a(z_a) = z_a + c$. In our study, $\Sigma$ partitions into
two subsets $\Sigma_0$ and $\Sigma_{\infty}$. Each point from
$\Sigma_0$ has a cone angle of $2 \pi N$, where $N$ is a positive
integer. Each point $p_{\infty}$ form $\Sigma_{\infty}$ has an
open neighborhood $W' \subset X$ with a map $\varphi_{\infty} : W'
\setminus \{p_{\infty}\} \to \CC$ such that $(W' \setminus
\{p_{\infty}\}, \varphi_{\infty})$ is a translation chart from the
atlas. Also, the set $\CC \setminus \varphi_{\infty}(W' \setminus
\{p_{\infty}\})$ is compact. Thus, the collection
$\Sigma_{\infty}$ contains all planar infinities on the surface.

Since all translations are holomorphic maps, the translation atlas
induces a complex structure on $X$ (for details see \cite{HS} and
\cite{M}). Moreover, the differential $dz_a$ in each $\varphi(W_a)
\subset \CC$ can be pulled back as a holomorphic differential
$\omega_a = \varphi_a^* dz_a$ in the corresponding $W_a$. But if
$$z_b = \varphi_b^{-1} \circ \varphi_a(z_a) = z_a + c$$ then
$dz_b=dz_a$. Hence, $\omega_a = \omega_b$ in any intersection $W_a
\cap W_b \neq \varnothing$ which gives rise to a global
holomorphic differential $\omega$ on $X \setminus \Sigma$.
Moreover, $\omega$ extends to the singular set $\Sigma$ so that
$\Sigma_0$ becomes the set of zeroes of $\omega$ and
$\Sigma_{\infty}$ becomes the set of all poles of $\omega$. The
latter are all double and with residue $0$. So we see that a
translation surface with planar infinities induces a pair
$(X,\omega)$ of a compact Riemann surface without boundary
together with an appropriate meromorphic differential.

To recover the translation atlas from a pair $(X,\omega)$, one can
cover $X \setminus (\text{zeroes}(\omega))$ with topological discs
$W_a$. On each of them define the chart $\varphi_a(p) =
\int_{p_a}^{p} \omega$, where $p_a \in W_a$ is fixed and $p$
varies in $W_a$. As $\omega$ is either holomorphic or meromorphic
with a double pole and residue $0$ inside the topological disc
$W_a$, the path of integration in $W_a \setminus
\text{poles}(\omega)$ is arbitrary. If $W_a \cap W_b \neq
\varnothing$ then $$z_b = \int_{p_b}^{p} \omega = \int_{p_a}^{p}
\omega + \int_{p_b}^{p_a} \omega = z_a + c$$ for $p \in W_a \cap
W_b $. Thus, we have obtained the desired translation atlas. As we
can see, the description of a translation surface with planar
infinities which we gave in the beginning of the current section
is equivalent to definition \ref{Definition_surface}.

The horizontal foliation $\Foliation_{\omega}$ on $X$, mentioned
in the introduction, is defined as follows. Let $\Foliation_{\CC}$
be the foliation of horizontal lines $\{z \in \CC \, |
\,\text{Im}(z) = s\}, \,\, s \in \Realline$ in $\CC$ oriented from
left to right (see figure \ref{Sector}b). Define the pulled-back
local foliation $\Foliation_a = \varphi_a^* \Foliation_{\CC}$ in
each $W_a$. Observe that $\Foliation_{\CC}$ is invariant with
respect to any translation, i.e. the translations map any
horizontal line to a horizontal line. Hence, $\Foliation_a =
\Foliation_b$ on each $W_a \cap W_b \neq \varnothing$. Thus, all
local foliations fit together in a global foliation
$\Foliation_{\omega}$ on $X$ with geodesic leaves and
singularities $\Sigma$. The oriented leaves of
$\Foliation_{\omega}$ determine globally a horizontal direction on
$(X,\omega)$. Since translations are Euclidean isometries, the
Euclidean metric on $\CC$ induces a Euclidean metric on $X
\setminus \Sigma$. In this metric geodesics that do not go through
singularities are isometric to straight lines in $\CC$. The notion
of a direction at a non-singular point $p \in X$ is as defined in
the introduction. It is the counterclockwise angle between the
horizontal leaf and an oriented geodesic both passing through $p$.
Finally, an oriented geodesic always forms the same angle with any
horizontal leaf it intersects, so it never self-intersects, except
possibly to close up.

\begin{figure}[h!]
\centering
\includegraphics[width=8cm]{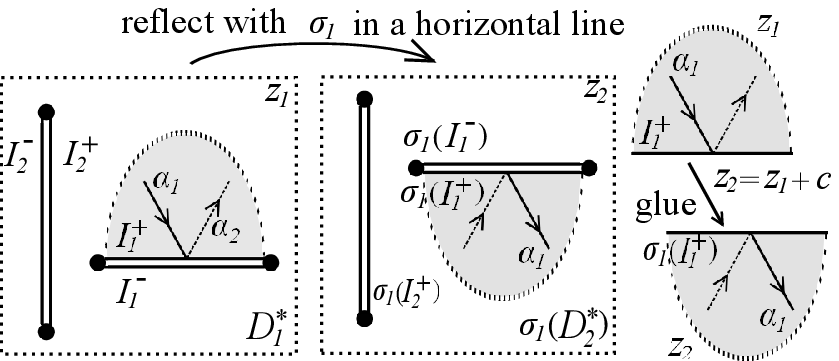}
\caption{} \label{Gluing}
\end{figure}

\paragraph{Construction.} Assume we have a configuration of disjoint compact line
segments $I_1,...,I_m$ in the plane $\CC$, which we regard as
two-sided mirrors. The angle between any two of them is a
rational-multiple of $\pi$. Observe that if one of the mirrors
forms a rational-$\pi$ angle with the rest of the mirrors, then
immediately follows that any pair of mirrors forms a
rational-$\pi$ angle. This is a consequence of the fact that in an
Euclidean triangle the angles at the vertices sum up to $\pi.$

To understand better the construction that follows, one could have
a simple toy-example in mind. Let us have two perpendicular
mirrors $I_1$ and $I_2$ in the plane $\CC$ like the ones depicted
on figure \ref{Gluing}. 

\begin{figure}[h!]
\centering
\includegraphics[width=8cm]{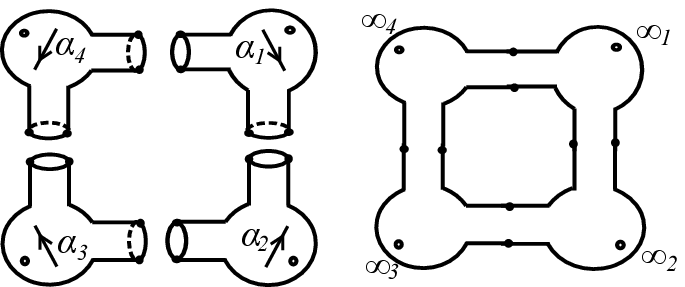}
\caption{} \label{SurfaceGluing}
\end{figure}

Begin by slicing $\CC$ along the segments $I_1,..., I_n$ to obtain
a closed slitted domain $D^*$ in which every mirror segment $I_k$
is doubled in order to obtain two parallel copies $I_k^+$ and
$I_k^-$ that form the boundary component of the surface $D^*$
around the slit $I_k$. For an intuitive geometric picture of $D^*$
in the case of the toy-example, look at figure \ref{Gluing}. Then
$D^*$ is homeomorphic to a once-punctured sphere with $n$ disjoint
open discs removed, as shown on figure \ref{SurfaceGluing} for the
case of two orthogonal mirrors. In particular, $\partial D^* =
\sqcup_{k=1}^n (I_k^+ \cup I_k^-)$.

For each segment $I_k$, fix the line $l_k \subset \CC$ through $0
\in \CC$ parallel to $I_k$. Denote by $\sigma_k$ the reflection of
$\CC$ in $l_k$. The group $G$ generated by all $\sigma_k, \,\, k
=1,..,n$ is a finite group. If $\alpha_1$ is a generic direction
in $\CC$, then $G(\alpha_1) =\{g(\alpha_1) \,|\, g \in G\} =
\{\alpha_1,..., \alpha_m\}$ is an orbit of maximal length $m \leq
n$. In our example $G \cong \Integers_4$ and a generic orbit has
$4$ elements. Pick $m$ copies $D^*_j$ of $D^*$ each with a choice
of a direction $\alpha_j$ in it. If you prefer more formally, let
$D_j^* = (D^*,\alpha_j)$. On figure \ref{SurfaceGluing}, in the
case of the toy-example, we can see a topological model of these
four slitted planes with a choice of direction on each of them. 
We glue $D_i^*$ to $D_j^*$ if and only if there is a segment $I_k
\subset \CC$ whose corresponding reflection $\sigma_k$ satisfies
$\sigma_k(\alpha_i)=\alpha_j$. The gluing is done in the following
way. Take $D_i^*$ and $\sigma_k(D_j^*)$. 
Glue the edge $I^+_k \subset D^*_i$ to the edge $\sigma_k(I^+_k)
\subset \sigma_k(D^*_j)$ and the edge $I^-_k \subset D^*_i$ to the
edge of $\sigma_k(I^-_k) \subset \sigma_k(D^*_j)$. On figure
\ref{Gluing} of the toy-example, we have chosen $i=1$ and $j=2$.
The upper edge $I^+_1 \subset D_1$ of the cut $I_1$ is glued to
the lower edge $\sigma_1(I_1^+) \subset \sigma_1(D^*_2)$ of the
cut $\sigma_1(I_1)$. Analogously, the lower edge $I^-_1$ from
$D^*_1$ is glued to upper edge $\sigma(I^-_1)$ from
$\sigma_1(D^*_2)$.

Both $D_i^*$ and $\sigma_k(D_j^*) $ are naturally translation
surfaces with piecewise geodesic boundaries, global coordinates
$z_i$ and $z_j$, and differentials $dz_i$ and $dz_j$ respectively.
Segments $I_k$ and $\sigma_k(I_k)$ are equal and parallel, hence
the gluing map is a translation $z_j = z_i + c$ (see the gluing of
the shaded pieces on figure \ref{Gluing}). Therefore the resulting
surface made out of $D^*_i$ and $\sigma_k(D^*_j)$ has a
translation structure. Moreover, $dz_j = dz_i$ along the gluing
locus, so there is a well-defined holomorphic differential on the
new surface which extends meromorphically to both of its infinity
points.

Now, follow the described gluing procedure for all cuts on the
pieces $D^*_j$, where $j=1,..,m$. The final result is a closed
Riemann surface $X$ and a meromorphic differential $\omega$ with
only double poles and zero residues, as well as simple zeroes with
cone angle $4 \pi$. For the example of the two orthogonal mirrors,
figure \ref{SurfaceGluing} illustrates how the four pieces
$D_1^*,...,D_4^*$ fit together to form a compact torus $X$ with a
complex structure and a meromorphic differential $\omega$ on $X$.
There are eight simple zeroes of $\omega$ and four double poles.
The zeroes are obtained from identifying pairs of black vertices
on the segments $I_k$ form figure \ref{Gluing}. The cone angle at
each zero is $4\pi$ and the residue at each pole is $0$ as
desired.

\section{Proofs}
\paragraph{Proof of theorem \ref{Main theorem about directions}.}
From now on $(X,\omega)$ is an arbitrary translation surface with
planar infinities and $p_0 \in X \setminus (\text{zeroes}(\omega)
\cup \text{poles}(\omega)$ any fixed point. The idea is to cut out
a rectangle around each pole $\infty_j \in \text{poles}(\omega)$
and replace it by a one-handle.
\begin{figure}
\centering
\includegraphics[width=8cm]{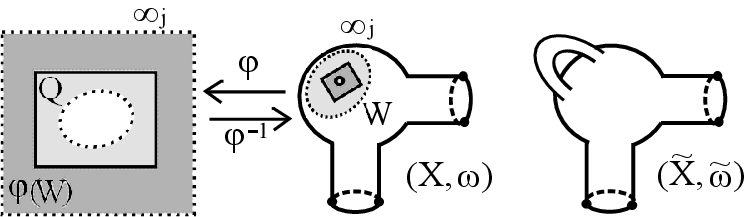}
\caption{} \label{Surgery}
\end{figure}
Indeed, choose a small topological disc $W$ around $\infty_j$ and
map it to $\CC$ by $\varphi(p) = \int_{q_0}^{p} \omega$ where $p$
varies in $W$ and $q_0 \in W$ is fixed. Notice, $\varphi$ is well
defined as the residue at $\infty_j$ is $0$, so the path of
integration is irrelevant. The image $\varphi(W) \subset \CC$ is
the complement of a compact set (the total shaded region on figure
\ref{Surgery} stretching to infinity). Draw a rectangle $Q \subset
\varphi(W)$ as shown on figure \ref{Surgery} and remove its
exterior (the darker region). On the surface, we remove the darker
rectangular domain containing $\infty_j$. Then glue together the
lower horizontal edge of $Q$ to the upper and the left to the
right, like gluing a torus. The gluing maps are clearly a vertical
and a horizontal translation respectively. Therefore we obtain a
handle with a translation structure compatible with the structure
on the rest of the surface (see figure \ref{Surgery}). By doing
this for each $\infty_j$, we obtain a compact translation surface
$(\tilde{X}, \tilde{\omega})$ of $\text{genus}(\tilde{X}) =
\text{genus}(X) + \sharp(\text{poles}(\omega))$, where
$\tilde{\omega}$ is now holomorphic (has no poles). A lot is known
about the behavior of the geodesics on such surfaces \cite{HS},
\cite{M}, \cite{V}, so we use this knowledge in our advantage. Let
$\tilde{\Lambda}_{p_0}$ be the set of all directions $\theta \in
\Circle$ for which the geodesic $\tilde{l}_{p_0}(\theta)$ on
$\tilde{X}$ is closed or hits a zero of $\tilde{\omega}$. Also,
let $\tilde{\Xi}$ be the set of all directions $\theta \in
\Circle$ for which the geodesic flow of $(\tilde{X},
\tilde{\omega})$ in direction of $\theta$ is minimal \cite{M}
(e.g. an ergodic flow is minimal \cite{HS},\cite{M}). Then
$\tilde{\Lambda}_{p_0}$ is countable but dense in $\Circle$ (see
\cite{V}) and $\tilde{\Xi}$ is dense and of full measure in
$\Circle$ (see \cite{M}, \cite{HS}). As a result, the set
$\tilde{\Theta}_{p_0} = \tilde{\Xi} \setminus
\tilde{\Lambda}_{p_0}$ consists of all $\theta \in \Circle$ for
which the geodesic ray $\tilde{l}_{p_0}(\theta)$ is dense in
$\tilde{X}$. Moreover, $\tilde{\Theta}_{p_0}$ is dense and of full
measure in $\Circle$. Therefore, for any $\theta \in
\tilde{\Theta}_{p_0}$ the corresponding geodesic ray
$l_{p_0}(\theta)$ on the original surface $(X,\omega)$ hits a pole
of $\omega$.

Let $U_{p_0} \subset \Circle$ be the set of all directions $\theta
\in \Circle$ with the property that the geodesic ray
$l_{p_0}(\theta)$ on $(X,\omega)$ in the direction of $\theta$
reaches a pole of $\omega$. Since the geodesic flow on
$(X,\omega)$ depends continuously on the initial point and
direction, the condition that a geodesic ray reaches a planar
infinity is open. Therefore, for each $\theta \in U_{p_0}$ there
exists an open circular interval $(\alpha,\beta) \subset U_{p_0}$
that contains $\theta$ and for any $\theta' \in (\alpha,\beta)$
the ray $l_{p_0}(\theta')$ also reaches the same infinity. Hence,
$U_{p_0}$ is open in $\Circle$. Moreover, the dense set of full
measure $\tilde{\Theta}_{p_0}$ is contained in $U_{p_0}$.
Therefore, $U_{p_0}$ is open and dense set of full measure in
$\Circle$.

The second part of theorem \ref{Main theorem about directions}
follows from the first one. If we are given a rational mirror
configuration, unfold it into a translation surface with planar
infinities $(X,\omega)$ as described earlier. Then, the infinity
of the mirror domain lifts to the set of poles of $\omega$ on $X$
and we apply the first part of the theorem.

\begin{figure}
\centering
\includegraphics[width=8cm]{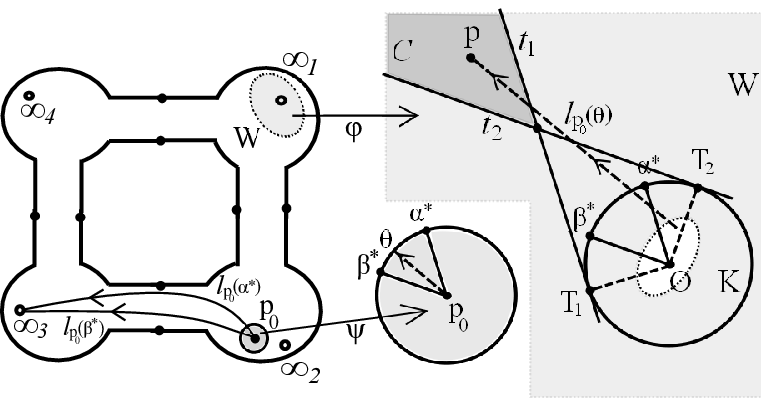}
\caption{} \label{ProofByPicture}
\end{figure}

\paragraph{Proof of theorem \ref{Main theorem illumination on surface}.}
As an open dense subset of $\Circle$, the constructed $U_{p_0}$ is
a countable disjoint union of open circular intervals $(\alpha_j,
\beta_j) \subset \Circle$, i.e. $U_{p_0} = \sqcup_{j=1}^{\infty}
(\alpha_j, \beta_j)$. By construction, the geodesic rays
$l_{p_0}(\theta)$ emitted from $p_0$ in all directions $\theta \in
(\alpha_j, \beta_j)$ go to the same pole of $\omega$. Fix some $j$
and take a subinterval $(\alpha^*,\beta^*) \subseteq (\alpha_j,
\beta_j)$ (it may even be convenient to choose $(\alpha^*,\beta^*)
= (\alpha_j, \beta_j)$). Choose $(\alpha^*,\beta^*)$ so that its
measure is less than $\pi$. Notice, that for every $\theta \in
(\alpha^*,\beta^*)$, each ray $l_{p_0}(\theta)$ on $X$ goes to the
same $\infty^* \in \text{poles}(\omega)$. In particular, $\infty^*
= \infty_3$ on figure \ref{ProofByPicture}. As
$\sharp(\text{poles}(\omega)) \geq 2$, take another $\infty \in
\text{poles}(\omega) \setminus \{\infty^*\}$ and call it
$\infty_1$ just like on our picture below. Choose a "small"
topological disc $W$ around $\infty_1$ with the property $W \cap
(\text{zeroes}(\omega) \cup \text{poles}(\omega)) = \{\infty_1\}$.
Define the translation chart $\varphi(p) = \int_{q_0}^{p} \omega$,
where $p$ varies in $W$ and $q_0 \in W$ is fixed. The zero residue
at $\infty_1$ guaranties independence of the integral on the path
between $q_0$ and $p$ in $W$. On figure \ref{ProofByPicture} we
have also provided an analogous chart $\psi$ around $p_0$. From
now on, we use the same notations in $W$ as the ones in
$\varphi(W)$. Thus, we identify $W$ with $\varphi(W)$. In $\CC$
the domain $W$ looks like the complement of a compact set (the
shaded region on figure \ref{ProofByPicture}). Let $K \subset W$
be a Euclidean circle in $\CC$ centered at $O$ and containing $\CC
\setminus W$ in its interior.

Abusing notation, let $\alpha^*$ and $\beta^*$ be the two points
on the circle $K$ such that the counter-clockwise angles between
the positive horizontal line through $O$ in $\CC$ and the radii
$O\alpha^*$ and $O\beta^*$ are respectively $\alpha^*$ and
$\beta^*$. Let points $T_1$ and $T_2$ on $K$ be such that
counter-clockwise $\measuredangle \alpha^*OT_1 = \measuredangle
T_2O\beta^* = \frac{\pi}{2}$. Draw the lines $t_1$ and $t_2$
tangent to circle $K$ at $T_1$ and $T_2$ respectively. Then they
bound an infinite sector $C$, depicted on figure
\ref{ProofByPicture} as a darker shaded region.

We claim that that  $C \subset X$ is not illuminated by $p_0$.
Assume that for some point $p \in C$ there exists $\theta \in
\Circle$ such that the geodesic $l_{p_0}(\theta) \subset X$
staring from $p_0$ in the direction of $\theta$ passes through
$p$. Then, clearly $l_{p_0}(\theta)$ goes to $\infty_1$. As
already commented in the introduction, any smooth geodesic on a
translation surface forms the same angle with the horizontal
direction at every point it passes through. In particular, the
angle between $l_{p_0}(\theta)$ and the horizontal direction in
the chart $W$ as well as near the point $p_0$ is always $\theta$.
By looking at the picture of the chart $W$ on figure
\ref{ProofByPicture}, we see that $\theta \in (\alpha^*,\beta^*)$
in $W$. Hence $\theta \in (\alpha^*,\beta^*) \subset \Circle$ at
the point $p_0$ as well. By the choice of the circular interval
$(\alpha^*,\beta^*)$, the geodesic ray $l_{p_0}(\theta)$ should go
to $\infty^* \neq \infty_1$. But a geodesic ray can only reach one
pole of $\omega$, so we get to a contradiction. Therefore, the
infinite sector $C$ on $(X,\omega)$ is not illuminated by $p_0 \in
X$.

To conclude the proof, notice that for each circular interval
$(\alpha^*,\beta^*) \subset U_{p_0}$ the unilluminated sector $C$
near $\infty_1$ can be also constructed around any other pole
$\infty \neq \infty^*$ of $\omega$, i.e. there are $k-1$
unilluminated copies of $C$. Partition $U_{p_0}$ into disjoint
subintervals for which we can apply the construction of
unilluminated infinite sectors from the preceding two paragraphs.
Thus, the the total sum of the angles of all unilluminated sectors
constructed on $(X,\omega)$ is $k-1$ times the total measure of
$U_{p_0} \subset \Circle$ which is $2 \pi$. Hence, the total angle
is $2 \pi(k-1)$.

\paragraph{Proof of theorem \ref{Main theorem illumination in
plane}.} Let $D \subset \CC$ be a rational mirror domain and $p_0
\in D$ (see figure \ref{Mirrors} or \ref{PictureCircleMap}).
Recall the finite group $G$ generated by all reflections in the
lines through $0\in \CC$ parallel to the mirrors. It acts on
$\Circle$ by rotations. Let $f_{p_0} : U_{p_0} \to \Circle$ be the
map described at the end of subsection "Main results" (see also
figure \ref{PictureCircleMap}) and assume it is not injective.
Then, there are $\theta_1 \neq \theta_2$ from $U_{p_0}$ such that
$f_{p_0}(\theta_1) = f_{p_0}(\theta_2)$. Take the finite orbit
$G(\theta_1) = \{g(\theta) \in \Circle \,\,| \,\, g \in G\}$. Then
$\theta \in G(\theta_1)$ if and only if $f_{p_0}(\theta) \in
G(\theta_1)$ so $\theta_2 \in G(\theta_1)$. Hence, the restriction
$f_{|_{G(\theta_1)}} : G(\theta_1) \to G(\theta_1)$ is not
bijective and there is $\theta^* \in G(\theta_1)$ such that
$\theta^* \in U_{p_0} \setminus f_{p_0}(U_{p_0})$. Since $f_{p_0}$
is a restriction of a rotation on each connected component of
$U_{p_0}$, there is $(\alpha^*,\beta^*) \ni \theta^*$ such that
$(\alpha^*,\beta^*) \subset U_{p_0} \setminus f_{p_0}(U_{p_0})$.
Remember the circle $K$ from figure \ref{PictureCircleMap} that
encompasses the mirrors and $p_0$. Using the circular interval
$(\alpha^*, \beta^*)$, we can carry out absolutely the same
construction as the one in the chart $W$ described in the proof of
theorem \ref{Main theorem illumination on surface}. For a picture
of this construction look at the rightmost large shaded area $W$
on figure \ref{ProofByPicture}. Observe that the notations of the
current proof match the picture's notations so that we can use it
directly, thinking that the set of mirrors is in the little white
elliptic region containing the center $O$. We claim that the
infinite sector $C$ (the darker shaded area) is not illuminated by
the source $p_0 \in D$. Indeed, assume there is a light ray
emitted by $p_0$ that reaches some $p \in C$. Then, from the
picture, the direction of this ray is $\theta \in
(\alpha^*,\beta^*)$. But the light ray started from $p_0$ in some
direction $\theta_0 \in \Circle$, so $\theta = f_{p_0}(\theta_0)$
which is a contradiction.

\end{document}